\documentclass[11pt]{amsart}
\def\Im{\textrm{Im}}

\def\Re{\textrm{Re}} 
\newtheorem{theoreme}{Theorem}

\newtheorem{proposition}{Proposition}
\newtheorem{lemme}[proposition]{Lemma}

\newtheorem{definition}[proposition]{Definition}
\newtheorem{remarque}[proposition]{Remark}

\numberwithin{equation}{section}
\numberwithin{proposition}{section}

\begin{document}
\title[Smoothing Effect for Schr\"odinger Equations]{Smoothing effect for Schr\"odinger\\ boundary value problems}

\author{N. Burq}
\address{Universit{\'e} Paris Sud,
Math{\'e}matiques,
B{\^a}t 425, 91405 Orsay Cedex}
\email{Nicolas.burq@math.u-psud.fr}
\urladdr{http://www.math.u-psud.fr/$\tilde{\ }$burq}

%\date{October $8^{\text{th} }$, 2002}

%\dedicatory{Preliminary version, revised \today}
%\keywords{Differential geometry, algebraic geometry}

\begin{abstract}We show the necessity of the non trapping condition for the plain smoothing effect ($H^{1/2}$) for Schr\"odinger equation with Dirichlet boundary conditions in exterior problems. We also give a class of trapped obstacles (Ikawa's example) for which we can prove a weak ($H^{1/2 - \varepsilon}$) smoothing effect\\
{\sc R\'esum\'e.} On d\'emontre que l'hypoth\`ese de non capture est  n\'ecessaire  pour l'effet r\'egularisant ($H^{1/2}$) pour l'\'equation de Schr\"odinger avec conditions aux limites de Dirichlet \`a l'ext\'erieur d'un domaine de ${\mathbb R}^d$. On donne aussi une classe d'obstacles captifs (l'exemple d'Ikawa) pour lesquels on d\'emontre un effet r\'egularisant affaibli ($H^{1/2 - \varepsilon}$).  
\end{abstract}
\maketitle
\section{Introduction}
Consider $u=e^{it\Delta} u_0$ solution of the
Schr\"odinger equation
\begin{equation}\left\{
\begin{aligned}
(i \partial _t+\Delta) u&=0 \text{ in $\mathbb{ R}\times\mathbb{ R}^d$},\\
u|_{t=0}&=u_0\in L^2({\mathbb{R}}^d) .
\end{aligned}\right.
\end{equation}
It is well known that $u\in L^\infty ({\mathbb{R}}_t;L^2({\mathbb{R}}^d))$
satisfies the following smoothing effect (for any  $s>1/2$ if $d\geq 3$)
\begin{equation}\label{smot}
  \Vert u\Vert _{L^2({\mathbb{R}}_t;\dot H^{1/2}_s({\mathbb{R}}^d))}\leq C \Vert u_0\Vert _{L^2},
\end{equation} 
where 
\begin{equation}
\dot H^{1/2}_s=\{u\in{\mathcal D}'({\mathbb{R}}^d);\langle x\rangle^{-s}\Delta^{1/4}
u\in L^2({\mathbb{R}}^d)\}.
\end{equation}
This result, which can be proved by explicit calculations, has been extended to more complicated operators, satisfying a non trapping assumption (see the results of Constantin - Saut~\cite{CoSa89}, Ben-Artzi - Devinatz~\cite{BADe98},  Ben-Artzi - Klainerman~\cite{BAKl92}, Do\"\i~\cite{Do96,Do96-1}, and Kato - Yajima~\cite{KaYa89}). It has been recently extended to the case of boundary value problems by G\'erard, Tzvetkov and the author~\cite{BuGeTz02}.\par

On the other hand, in~\cite{Do00} Do{\"\i}  has proved that, for Schr\"odinger operators in ${\mathbb R}^d$, the non trapping assumption is necessary for the $H^{1/2}$ smoothing effect. \par
In this paper we extend this latter result to the case of boundary value problems. Our first result reads as follows:
\begin{theoreme} 
\label{th2}Consider an arbitrary smooth domain with boundary  $\Omega\subset {\mathbb R}^d $, with no infinite order contact with its boundary (see the precise definition in Section~\ref{sec.4}),  and $P$ a  second order self adjoint operator on $L^2(\Omega)$, with domain $D\subset  H^1_{0}( \Omega)$ and such that the boundary is non characteristic.
Denote by $\varphi_{s}: \ ^bT^*\Omega\setminus\{0\}\rightarrow \ ^bT^*\Omega\setminus\{0\}$  the bicharacteristic flow of the operator $P$ (given by the integral curves of the Hamiltonian vector field of  the principal symbol of $P$ reflecting on the boundary according to the law of geometric optics--see Section~\ref{sec.4}--) defined on the boundary cotangent bundle. Let $A\in \Psi^{(1/2)}$ be a classical tangential pseudodifferential operator of order $1/2$.  Suppose that $(z_{0}, \zeta_{0})\in \ ^ bT^* \Omega\setminus\{0\}$ satisfy the trapping assumption:
\begin{equation}
\label{trap}\int_{-\infty}^{0}|\sigma_{1/2}(A)(\varphi_{s}( x_{0}, \zeta_{0})| ^2 ds= + \infty
,\end{equation}
where $\sigma_{1/2}(A)$ is the principal symbol of the operator $A$.
Then for any $t_{0}>0$ the map  
\begin{equation}
\label{eq1.10}u_{0}\in C^\infty_{0} \subset L^2 (\Omega)\mapsto A e^{itP}u_{0}\in L^{2}([0, t_{0}];L^2(M))
\end{equation} is not bounded (even for data with fixed compact support).    
\end{theoreme}

\begin{remarque}The assumption~\eqref{trap} can be essentially fulfilled in two distinct cases
\begin{enumerate}
\item If $A$ is compactly supported (in the $z$ variable), then~\eqref{trap} means that 
the bicharacteristic starting from $(z_{0}, \zeta_{0})$ spends an infinite time in the support of $A$, which corresponds to a ``trapped trajectory''
\item If $A$ is not compactly supported, a typical example is (in the case $P= - \Delta$) $A(z, D_{z})= a(|z|) |D_{z}|^{1/2}$, then~\eqref{trap} might correspond to a lack of decay of $a(x)$ at infinity: suppose that the trajectory starting from $(z_{0}, \zeta_{0})$ is not trapped; hence it leaves any compact set and for $\pm s \rightarrow + \infty $, $(z(s), \zeta(s))\sim (s\zeta_{\pm}, \zeta_{\pm})$ and~\eqref{trap} is equivalent to $|a|^2\notin L^1( {\mathbb R})$ (and we recover the usual assumption required for proving the smoothing effect, see~\cite{Do96-1}).    
\end{enumerate}
\end{remarque}
\begin{remarque}
  We could have added lower order terms to $P$ and supposed that the Cauchy problem is well posed in $L^2$ (in case of first order terms). The condition~\eqref{trap} has in this case to be modified.
\end{remarque}
\begin{remarque}
\label{re1}In~\cite{Do96-1, Do00}, Do{\"\i}  proves this result in the case of a manifold without boundary and gives some variants of this result for operators of higher order, and with weights in times. The proof we present below is essentially self contained in this case and it can also handle these variants modulo slight modifications. The proof in presence of a boundary is much more technical.
\end{remarque}
\begin{remarque}
For $P= -\Delta_{g}$ the $x$-projection of the integral curves of $H_{p}$ are the geodesics for the metric $g$. 
\end{remarque}
\begin{remarque}
The smoothness assumption can be relaxed to $C^2$ coefficients and $C^3$ domains (and even to $C^1$ coefficients, but the assumption~\eqref{trap} is then more complicated since the Hamiltonian flow is no more well defined). We also can prove Theorem~\ref{th2} for systems (see Remark~\ref{re4.20}).
\end{remarque}
Having Theorem~\ref{th2} in mind, a natural question is  whether a weakened version of~\eqref{smot} might hold for some trapping geometries. In the case of a stable (elliptic) trapped trajectory, the existence of quasi-modes well localized along this trajectory shows that no such result may hold (see Remark~\ref{re4.7}). However in the case of hyperbolic trapped trajectories, we do obtain such a weak smoothing effect:
\begin{theoreme}
\label{th6}Consider $\Theta= \cup _{i=1}^N \Theta_{i} \subset {\mathbb R}^d$ a finite union of strictly convex obstacles satisfying the assumptions of Section~\ref{sec.2}. Denote by $\Omega= \Theta^c$ its complement. Then  for any $\varepsilon>0$ and $\chi \in C^\infty _{0}({\mathbb R}^d)$ there exists $C>0$ such that for any $u_{0}\in L^2( \Omega)$,
\begin{equation}\label{weaksmooth}
\| \chi e^ {it \Delta_{D}}u_{0}\|_{L^2( {\mathbb {R}}_{t}; H^{1/2 -\varepsilon}( \Omega)}\leq C \|u_{0}\|_{L^2( \Omega)}
\end{equation}
\end{theoreme}
\begin{remarque}
This result was proved in~\cite{BuGeTz02} with no $\varepsilon$ loss under the non trapping assumption: ``any geodesic of the metric $g$ reflecting on the boundary according to the laws of geometric optics goes to the infinity'' which is clearly no fulfilled here.
\end{remarque}  
To prove Theorem~\ref{th2} we will follow the same kind of strategy as in Do\"\i 's paper~\cite{Do96-1, Do00}. However we will replace in his argument the use of Egorov's theorem by the use of the theorem of propagation of Wigner measures, which has three advantages: first it simplifies the rest of the proof, second it allows to relax assumptions (on the regularity of the coefficients) and finally the proof holds also for a (system of) boundary value problem (whereas Egorov's Theorem is not true in these cases)\par
To prove Theorem~\ref{th6} we reduce, following~\cite{BuGeTz02}, the estimate~\eqref{weaksmooth} to obtaining estimates for the outgoing resolvent of $\Delta_{D}$, $(- \Delta- (z\pm i0))^{-1}$. Then we show that these estimates can be deduced from a combination of other estimates proved by M. Ikawa~\cite{Ik82, Ik83,Ik88} and some form of the maximum principle.\par
  The article is written as follows:  in Section~\ref{sec.3} we recall the definition of Wigner measures which will be used in the sequel and we prove Theorem~\ref{th2} in the simpler case where $\Theta= \emptyset$. In Section~\ref{sec.4} we give the necessary modifications required to handle the general case $\Theta\ne \emptyset$. In Section~\ref{sec.2} we prove Theorem~\ref{th6}. Finally we have stated at the end of Section~\ref{sec.2} an application of our smoothing result to the global existence of non linear Schr\"odinger equations.
\par
\noindent{ \bf Acknowledgements :} I would like to thank C. Zuily for discussions about the results in this article and the referees whose observations lead to substantial improvements in the exposition. {This work was completed during a stay  at the Department of Mathematics of the University of California, Berkeley, partially funded by a N.A.T.O. fellowship and the Miller Institute for Basic Research in Sciences. I thank these institutions.} 
\section{Proof of Theorem~{\protect\ref{th2}}: the case of empty boundary}\label{sec.3}
\subsection{Wigner measures}
In this section we recall the definition of Wigner measures (or semi-classical measures) introduced by G{\'e}rard -- Leichtnam~\cite{GeLe93-1} and  Lions -- Paul~\cite{LiPa93} (see also the survey by G{\'e}rard -- Markowich -- Mauser -- Poupaud~\cite{GMMP97}). We work in the context of functions of $1+d$ variables ($(t,z)$) in  $ L^2_{\text{loc}}( {\mathbb R}_{t}; L^2( {\mathbb R}^d_{z}))={\mathcal {L}}^2$ and we have adapted the definitions in~\cite{GeLe93-1, LiPa93} to fit our purpose. 
\begin{definition}
We will say that a sequence of functions $(f_{n} )\in {\mathcal {L}}^2$ is bounded in ${\mathcal {L}}^2$ if for any $\varphi\in C^\infty _{0}( {\mathbb R}_{t})$, the sequence $(\varphi f_{n})$ is bounded in $L^2$.
\end{definition}
\begin{definition}
 We will say that an operator $A$ is bounded on ${\mathcal {L}}^2$ if there exists $\varphi\in C^\infty_{0}( {\mathbb R}_{t})$ such that for any $f\in {\mathcal {L}}^2$,
 $$\| Af \|_{L^2_{t,z}}\leq C \| \varphi f \| _{L^2_{t,z}}.$$
\end{definition}
Denote by $(x, \xi) = (t,z,\tau, \zeta)$ a point in $T^*{\mathbb R}^{d+1}$; and  
consider for  $a(x, \xi)\in C^\infty_{0}( {\mathbb R}^{ 2d+2})$ and $\varphi\in C^\infty_{0}( {\mathbb R}_{t})$ equal to $1$ near the support of $a$,  the operator $\text{ Op}_{ \varphi}(a)(x, hD_{x})$ defined on ${\mathcal {L}}^2$ by
\begin{equation}
\label{eq2.6bis}
\begin{aligned}
\text{Op}_{\varphi}(a)(x, hD_{x})f&=\text{Op}_{\varphi}(a)(t,z,hD_{t}, hD_{z})f\\
&= \frac{ 1}{ (2 \pi)^{d+1} }\int e^{i(t\cdot \tau+ z\cdot \zeta)}a(t,z,h\tau, h\zeta) \widehat { \varphi(t) f}(\tau, \zeta)d\tau d\zeta.
\end{aligned}
\end{equation}
The operator $\text{ Op}(a)_{\varphi}(t,z, hD_{t}, hD_{z})$ is (uniformly with respect to $0<h<1$) bounded on ${\mathcal {L}}^2$  and we have the following weak form of the G{\aa}rding inequality
\begin{proposition}
\label{prop2.0}
For any $a\in C^\infty_{0}( {\mathbb R}^{ 2d+2})$ and any sequence $(f_{n})$ bounded in ${\mathcal {L}}^2$ and $(h_{n})\in ]0,1]$; $\lim_{n\rightarrow + \infty}h_{n}=0$, 
\begin{equation}
\label{eq2.7}a(x, \xi)\geq 0 \Rightarrow \liminf_{n\rightarrow + \infty}\Re \left(\text{ Op}(a)_{\varphi}(x,h_{n}D_{x})f_{n}, f_{ n}\right)_{L^2( {\mathbb R}^{d+1})}\geq 0.
\end{equation}
\end{proposition}
To prove this result consider for $\varepsilon>0$, $\psi\in  C^\infty_{0}( {\mathbb R}^{ 2d+2})$ equal to $1$ near the $(t,z,\tau, \zeta)$ projection of the support of $a$ and $b= \varphi(t) \sqrt{ \varepsilon+a}\psi(t,z, \tau, \zeta)\in  C^\infty_{0}( {\mathbb R}^{ 2d+2})$. Then the symbolic calculus shows
\begin{multline}
0\leq\left(\text{ Op}(b)_{\varphi}^*\text{ Op}(b)_{\varphi}f_{n}, f_{n}\right)\\
= \left(\text{ Op}(a)_{\varphi}(x,h_{n}D_{x})f_{n}, f_{ n}\right)_{L^2}+ \varepsilon\left(\varphi(t) \psi^2(x,h_{n}D_{x}) \varphi(t)f_{n}, f_{n}\right)_{ L^2}+ {\mathcal {O}}( h_{n}),
\end{multline}
hence, taking the $\liminf$ and using that $\liminf (\alpha_{n}+ \beta_{n})\leq \liminf (\alpha_{n})+ \limsup (\beta_{n})$, we get  
\begin{equation}
\liminf_{n\rightarrow + \infty}\Re \left(\text{ Op}(a)(x,h_{n}D_{x})f_{n}, f_{ n}\right)_{L^2( {\mathbb R}^d)}+ \varepsilon \limsup _{n\rightarrow + \infty}\| \psi(x,h_{n}D_{x})\varphi(t) f_{n}\|^2 \geq 0.
\end{equation} 
When $\varepsilon>0$ tends to $0$ we obtain Proposition~\ref{prop2.0}.
\par
By the symbolic calculus, the operator $\text{ Op}(a)_{\varphi}$ is modulo operators bounded on ${\mathcal {L}}$ by ${\mathcal {O}}(h^\infty)$, independent of the choice of the function $\varphi$. For conciseness, we will drop in the sequel the index $\varphi$.
As in~\cite{GeLe93-1} (see also~\cite{Bu97-1}) we can prove:
\begin{proposition}
\label{prop2.1} Consider a sequence $(f_{n})$ bounded in ${\mathcal {L}}^2$. There exist a subsequence $(n_{k})$ and a positive Radon  measure on ${\mathbb R}^{ 2d+2}$, $\mu$, such that for any $a\in C^\infty_{0}( {\mathbb R}^{ 2d+2})$
\begin{equation}
\label{eq2.8bis}\lim_{k\rightarrow + \infty}\left(\text{Op}(a)(x, h_{n_{k}}D_{x})f_{n_{k}}, f_{n_{k}}\right)_{ L^2}= \langle \mu, a(x, \xi)\rangle.
\end{equation} 
\end{proposition}
The idea for extracting such a sequence is to fix $ a$ and consider the bounded sequence $(L(a)_{n})=\left(\text{Op}(a)(x, h_{n}D_{x})f_{n}, f_{n}\right)_{ L^2}$. By compactness we can extract a subsequence which converges. Iterating this process for a sequence $(a_{j})$ dense in $C^\infty_{0}$, we obtain, by diagonal extraction, a sequence $(f_{n_{k}})$ such that the limit exists for any $a_{j}$. By~\eqref{eq2.7} the limit defines a positive functional on a dense subset of $C^\infty_{0}$ (hence this limit is continuous for the $C^0$ topology). It is consequently a Radon measure and the limit~\eqref{eq2.8bis} exists for any $a\in C^\infty_{0}$. For the sake of conciseness we shall denote again by $(f_{n})$ the extracted subsequences.\par
The measure $\mu$ represents at points $(x_{0}, \xi_{0})$ the oscillations of the sequence $(f_{n})$ at point $x_{0}$ and scale $\xi_{0}/h_{n}$. The oscillations at frequencies smaller than $h_{n}^{ -1}$ are concentrated in $\{\xi_{0}=0\}$ whereas the oscillations at higher ($>> h_n^ {-1}$) frequencies are lost.
\subsection{Invariance of the Wigner measure}
\subsubsection{Elliptic regularity}\label{se.2.2.1}
Suppose that the sequence $(f_{n})$ is solution of the equation
\begin{equation}\label{eq2.200}
(ih_{n}\partial _{t} + h_{n}^ 2 P)f_{n}={\mathcal {O}}(h_{n})_{{\mathcal {L}}^2}
\end{equation}
Take $a\in C^\infty_{0}$ and consider first
\begin{equation}
\label{eq2.16}
 \left(\text{Op}(a)(x,h_{n}D_{x}) (ih_{n}\partial _{t} + h_{n}^ 2 P)f_{n}, f_{n}\right)_{ L^2}= o(1).
\end{equation}
Taking into account that the operator $\text{Op}(a)(x, h_{n}D_{x}) (ih_{n}\partial _{t} + h_{n}^ 2 P)$ is equal to $\text{Op}(a\times (- \tau + p(z, \zeta))(x, hD_{x})$ modulo an operator bounded by ${\mathcal {O}}(h_{n})$ on ${\mathcal {L}}^2$ and passing to the limit in~\eqref{eq2.16} we obtain 
\begin{equation}
\label{eq2.17}\langle \mu, a(x, \xi)(- \tau+ p(z, \zeta))\rangle=0
\end{equation}
from which we deduce:
\begin{proposition}
\label{elliptic}
The measure $\mu$ is supported in the semi-classical characteristic set of the operator:
\begin{equation}\label{eq2.17'}\text{Char}(ih_{n}\partial _{t} + h_{n}^ 2 P)= \{ (x, \xi)=(t,z,\tau, \zeta);\tau = p(z, \zeta)\}
\end{equation}
\end{proposition}
\begin{remarque}
\label{re2.10}Suppose that the sequence $(f_{n})$ is solution of the equation~\eqref{eq2.200}. Then for any $a\in C^\infty_{0}( {\mathbb R}^{2d})$, the function 
\begin{equation}
t \mapsto \left(\text{Op}(a)( z, h_{n} D_{z})f_{n}\mid_{t}, f_{n}\mid_{t}\right)_{L^2( {\mathbb R}^d_{z})}(t)
\end{equation}
is, according to~\eqref{eq2.200}, locally uniformly equicontinuous. Hence using Ascoli's theorem, it is possible to extract a subsequence $(f_{n_{k}})$ (independent of $t$) such that there exist a family of positive measures $\mu_{t}$ continuous with respect to $t$ and such that for any $t$ and any $a \in C^\infty_{0}({\mathbb R}^{2d})$  we have 
\begin{equation}
\lim _{n\rightarrow + \infty}\left(\text{Op}(a)( z, h_{n} D_{z})f_{n}, f_{n}\right)_{L^2( {\mathbb R}^d_{z})}(t)= \langle \mu_{t}, a\rangle.
\end{equation}
Of course, from $\mu_{t}$ one can recover the measure $\mu$ (assuming that the extracted sequences are the same):
\begin{equation}\label{eq3.13}
\mu = dt \otimes \delta_{\tau =  p(z, \zeta)}\otimes \mu_{t}.
\end{equation}   
\end{remarque}
\subsubsection{Propagation of the Wigner measure}
Suppose now that \begin{equation}
(ih_{n}\partial _{t} + h_{n}^ 2 P)f_{n}=o(h_{n})_{{\mathcal {L}}^2}.
\end{equation}
Consider the bracket ($P^*=P$)
\begin{equation}
\label{eq2.12bis}\begin{aligned}
h_{n}^{ -1}&\left([\text{Op}(a)(x,h_{n}D_{x}), ih_{n}\partial _{t} + h_{n}^ 2 P]f_{n}, f_{n}\right)_{ L^2}\\ 
&=h_{n}^{ -1}\left( (ih_{n}\partial _{t} + h_{n}^ 2P)\text{Op}(a)(x, h_{n} D_{x})f_{n}, f_{n}\right)_{ L^2}+ o(1),\\
&= o(1). 
\end{aligned}
\end{equation}
 Taking into account that the operator 
\begin{equation}
h_{n}^{ -1}[\text{Op}(a)(x,D_{x}), ih_{n}\partial _{t} + h_{n}^ 2 P]
\end{equation}
 is equal to 
\begin{equation}
\frac 1 i \{a, - \tau + p(z, \zeta)\}(x, h_{n}D_{x})+ {\mathcal {O}}(h_{n})_{{\mathcal {L}}( {\mathcal {L}}^2)},
\end{equation}
 where the Poisson bracket of $a$ and $q$, $\{ a, q\}$, is defined by
\begin{equation}\{ a, q\}= \nabla_{\tau,\zeta}a\cdot \nabla_{ t,z}q- \nabla_{t,z}a \nabla _{\tau, \zeta}q,
\end{equation}
 we can pass to the limit in~\eqref{eq2.12bis} and obtain:
\begin{equation}
\label{eq2.14}\langle \mu,\{ a, -\tau +p(z, \zeta)\} \rangle=0,
\end{equation}
or equivalently (with $H_{\tau -p(z, \zeta)}$ the Hamiltonian vector field of $\tau -p$)
\begin{equation}
\label{eq2.15bis}H_{\tau - p(z,\zeta)}(\mu) = (\partial _{t}- H_{p})\mu=0.
\end{equation}
Gathering Proposition~\ref{elliptic} and~\eqref{eq2.15bis}, we have proved
\begin{proposition}
\label{propag}
The measure $\mu$ is  invariant  along the integral curves of the vector field $H_{\tau -p}$ drawn on the surface $\{\tau =p(z, \zeta)\}$. Equivalently, if we denote by $\varphi_{s}$ the Hamiltonian flow of the function $p(z, \zeta)$ on $T^* {\mathbb R}^2$ and if $\mu_{t}$ is as in~\eqref{eq3.13}, we have the equality for any $s \in {\mathbb R}$
\begin{equation}
\mu_{s}= \varphi_{s}^* ( \mu_{0}) \Leftrightarrow \langle \mu_{s}, a\circ \varphi_{s}\rangle= \langle \mu_{0}, a\rangle. \end{equation} 
\end{proposition}
\subsection{Proof of Theorem~{\protect\ref{th2}} in the case $\Omega= {\mathbb R}^d$}
Take $(z_{0}, \zeta_{0})$ satisfying the assumption~\eqref{trap} and consider $ \varphi\in C^\infty_{0}( {\mathbb R}^d)$ such that $\int |\varphi| ^2= 1$ and 
\begin{equation}\label{eq3.30}
u_{0,n}=n^{d/4}\varphi(n^{1/2}(z-z_{0}))e^{ in(z-z_{0})\cdot \zeta_{0}}.
\end{equation}
 Denote by $v_{n}= e^{itP}u_{0,n}$ the  corresponding solution of the Schr\"odinger equation. To prove Theorem~\ref{th2}, we are going to show:
\begin{equation}\label{eq3.1}
\forall \varepsilon>0 \lim_{n\rightarrow + \infty}\|A(z, D_{z}) v_{n}\|_{L^2([0, \varepsilon]\times {\mathbb R}^d)}= + \infty
\end{equation}
if $A \in S^{1/2}( {\mathbb R}^{2d})$ satisfies the assumptions of Theorem~\ref{th2}.\par
 For this we compute, with $h_{n}= 1/n$ and $\Psi \in C^\infty_{0}({\mathbb R}^{2d})$, $0\leq \Psi \leq 1$ equal to $1$ near $0$ and $\alpha>0$ fixed,
\begin{equation}\label{eq3.31}
\begin{aligned}\|A(z, D_{z}) v_{n}\|^2_{L^2([0, \varepsilon]\times {\mathbb R}^d)}&=\int_{0}^\varepsilon \left(A^*(z, D_{z})A(z, D_{z}) v_{n}, v_{n}\right)_{L^2( {\mathbb R}^d_{z})}dt\\
&\geq \int_{0}^\varepsilon \left(\Psi(\alpha z, \alpha h_{n}D_{z})A^*(z, D_{z})A(z, D_{z}) \Psi(\alpha z, \alpha h_{n}D_{z})v_{n}, v_{n}\right)_{L^2( {\mathbb R}^d_{z})}dt\\
& \qquad -C\\
&\geq \int_{0}^\varepsilon \left(h_{n}^{-1}b^*(z, h_{n}D_{z})b(z, h_{n}D_{z})v_{n}, v_{n}\right)_{L^2( {\mathbb R}^d_{z})}dt - C ,
\end{aligned}
\end{equation}
with $b(z, \zeta)= \sigma_{1/2}(A)(z, \zeta) \Psi (\alpha z, \alpha\zeta)$.\par
But, for any $T$, if $n$ is large enough 
\begin{equation}\label{eq3.32}
\int_{0}^\varepsilon \left(h_{n}^{-1}b^*(z, h_{n}D_{z})b(z, h_{n}D_{z})v_{n}, v_{n}\right)_{L^2( {\mathbb R}^d_{z})}dt\geq\int_{0}^{h_{n} T} \left(h_{n}^{-1}b^*(z, h_{n}D_{z})b(z, h_{n}D_{z})v_{n}, v_{n}\right)_{L^2( {\mathbb R}^d_{z})} dt.
\end{equation}
Denote by $u_{n}(s,z) = v(h_{n}s, z)$, the solution of the semi-classical Schr\"odinger equation
\begin{equation}
(i h_{n}\partial _{s}+ h_{n}^2 P) u_{n}=0
\end{equation}
we obtain for any $T>0$
\begin{equation}\label{eq3.33}
\int_{0}^\varepsilon \left(h_{n}^{-1}b^*(z, h_{n}D_{z})b(z, h_{n}D_{z})v_{n}, v_{n}\right)_{L^2( {\mathbb R}^d_{z})}dt\geq\int_{0}^{ T} \left(b^*(z, h_{n}D_{z})b(z, h_{n}D_{z})u_{n}, u_{n}\right)_{L^2( {\mathbb R}^d_{z})}ds. 
\end{equation}
According to~\eqref{eq3.30}, the Wigner measure, $\mu_{0}$, of the sequence $(u_{n}\mid_{t=0})$ is equal to 
$$\delta_{(z, \zeta)=(z_{0}, \zeta_{0})}.$$
 From Proposition~\ref{propag} and~\eqref{eq3.13}, we deduce that the Wigner measure, $\mu_{s}$, of $(u_{n}\mid_{t=s})$ is equal to $\delta_{(z,\zeta)=\varphi_{-s}(z_{0}, \zeta_{0})}$, where $\varphi_{s}$ is the flow of $H_{p}$.\par
 Hence (for $T$ fixed)
\begin{multline}
\label{eq3.34}
\lim_{n\rightarrow + \infty}\int_{0}^{ T} \left(b^*(z, h_{n}D_{z})b(z, h_{n}D_{z})u_{n}, u_{n}\right)_{L^2( {\mathbb R}^d)}dt\\
\begin{aligned}& = \int_{0}^ T \langle \mu_{s}, b\rangle ds\\
&= \int_{0}^T |b|^2(\varphi_{-s}(z_{0}, \zeta_{0})) ds\\&= \int_{0}^T|\sigma_{1/2}(A) (\varphi_{-s}(z_{0}, \zeta_{0}))| ^2 |\Psi ( \alpha \varphi_{-s}(z_{0}, \zeta_{0}))| ^2  ds.
\end{aligned}
\end{multline}
From~\eqref{eq3.31},~\eqref{eq3.32},~\eqref{eq3.33} and~\eqref{eq3.34} we deduce (if $\alpha$ is chosen small enough) that for any $T>0$ and with a fixed constant $C$ independent of $T$:
\begin{equation}
\liminf_{n\rightarrow + \infty}\|A(z, D_{z}) v_{n}\|^2_{L^2([0, \varepsilon]\times {\mathbb R}^d)}\geq \int_{0}^T|\sigma_{1/2}(A) (\varphi_{-s}(z_{0}, \zeta_{0}))| ^2 ds-C
\end{equation}
Letting $T$ tend to the infinity (and using the assumption~\eqref{trap}), we obtain~\eqref{eq3.1}.
\section{Proof of Theorem~\ref{th2} for a Dirichlet problem}\label{sec.4}
In this section we are going to give the outline of the proof of Theorem~\ref{th2} in the general case. In fact the proof is essentially the same as in the previous section. The differences are that we have to define Wigner measures for sequences bounded in $L^2_{\text{loc}}( {\mathbb R}_{t}; L^2( \Omega))$ and prove the elliptic (Proposition~\ref{elliptic}) and propagation (Proposition~\ref{propag}) results for these measures. Then we will construct a sequence of initial data whose Wigner measure is $\delta_{(z_{0}, \zeta_{0})}$ where $(z_{0}, \zeta_{0})$ satisfies the assumption~\eqref{trap} and the sequence of solutions of the Schr\"odinger equation with these initial data will prove the result. Fortunately, all these constructions have already been done (see the works by G\'erard -- Leichtnam~\cite{GeLe93-1}, Miller~\cite{Mi97, MI00}, Burq -- Lebeau~\cite{BuLe01} and Burq~\cite{Bu02}) in some slightly different settings. All that we have to do is to adapt these constructions to our framework and to glue the pieces together.\par
 For the sake of completeness, we are going to give an outline of the constructions. However, we insist on the fact that in this section, most of the material is taken from the works cited above.
\begin{remarque}
\label{re4.20}
For simplicity, we have restricted the study to the case of a scalar equation; however, following~\cite{BuLe01}, it would not be much more difficult to prove the result for systems.
\end{remarque}
 \subsection{Geometry}
Denote by $M= {\mathbb R}_{t}\times \Omega$, $x=(t,z) \in M$ and by ${^b}T M$  the bundle of rank $d+1$ whose sections are the vector fields tangent to $\partial M$, $^bT^* M$ the dual bundle (Melrose's compressed cotangent bundle) and  $j : T^* M\rightarrow {^b}T^* M$ the canonical map. In any coordinate system where $M=\{ x= (x_{n}>0, x')\}$), the bundle $^bT M$ is generated by the fields $\frac{\partial }{ \partial x'}$,
$x_{n} \frac{\partial }{ \partial x_{n}}$ and $j$ is defined by
\begin{equation} 
  j(x_{n},x',\xi_{n} ,\xi' )=(x_{n},x',v=x_{n}\xi_{n}, \xi' ) .
 \end{equation}
Denote by $\text{Car} \widetilde{P}$ the semi-classical characteristic manifold of $\widetilde{P}= ih \partial _{t}+ h^2P$ and $Z$ its projection
\begin{equation} \label{eq:5}
\text{Car} \widetilde{P}=\left\{ (x,\xi )= (t,z,\tau, \zeta)\in T^* {\mathbb R}^d\mid _{\overline M} ; p(x, \xi) =\tau\right\}, \qquad  
Z=j(\text{Car} \widetilde{P}).
 \end{equation}
The set $Z$ is a locally compact metric space.\par
Consider, near a point $x_{0}\in \partial M$ a geodesic system of coordinates for which $x_{0}= (0,0)$, $M= \{(x_{n},x')\in {\mathbb R}^+\times {\mathbb R}^{ d}\}$ and the operator $\widetilde{P}$ has the form  (near $x_{0}$)
\begin{equation}
\label{3.22} \widetilde{P}= -h^2 D_{x_{n}}^2  + R(x_{n}, x',hD_{x'})+ h Q(x, hD_{x}), 
\end{equation}
with $R$ a second order tangential operator and $Q$ a first order operator.\par
We recall now the usual decomposition of $T^*\partial M$ (in this coordinate system). Denote by $r(x_{n}, x',\xi')$ the semi-classical principal symbol of $R$ and $r_{0}= r\mid _{x_{n}=0}$. Then $T^*\partial M$ is the disjoint union of $\mathcal{ E}\cup \mathcal { G} \cup \mathcal{ H}$ with 
\begin{equation}
\label{3.24}\mathcal{ E}= \{ r_{0}<0\}, \mathcal{ G}=\{r_{0}=0\}, \mathcal{ H}= \{r_{0}>0\}.
\end{equation}
 Remark that $j$ gives a natural identification between $Z\mid_{\partial M}$ and ${\mathcal {H}}\cup {\mathcal {G}} \subset T^* \partial M$.
In $\mathcal{G}$ we distinguish between the {\em diffractive} points $\mathcal{ G}^{2,+}=\{ r_{0}=0, r_{1}= \partial _{x_{n}}r\mid _{x_{n}=0}>0\}$ and the {\em gliding} points $\mathcal{ G}^{-}=\{ r_{0}=0, r_{1}= \partial _{x_{n}}r\mid _{x_{n}=0}\leq 0\}$. We will make the assumption ($\Omega$ has no infinite order contact with its tangents) that for any $\varrho_{0}\in T^*\partial M$, there exists $N\in {\mathbb N}$ such that
$$H_{r_{0}}^N (r_{1})\neq 0$$\par
 The definition of the generalized bicharacteristic flow, $\varphi_{s}$ associated to the operator $P$ is essentially the definition given in~\cite{MeSj82}:
\begin{definition}
A generalized bicharacteristic curve $\gamma(s)$ is a continuous curve from an interval $I \subset {\mathbb R}$ to $ Z$ such that
\begin{enumerate}
\item if $s_{0}\in I$ and $\gamma(s_{0})\in T^* M$ then close to $s_{0}$, $\gamma$ is an integral curve of the Hamiltonian vector field $H_{\widetilde p}$
\item If $s_{0}\in I$ and $\gamma(s_{0})\in {\mathcal {H}}\cup {\mathcal {G}}^{2,+}$ then there exists $\varepsilon>0$ such that for $0<|s-s_{0}|<\varepsilon$,  $x_{n}(\gamma(s))>0$
\item If $s_{0}\in I$ and $\gamma(s_{0})\in{\mathcal {G}}^{-}$ then for any function $f \in C^\infty(T^*{\mathbb R}^{d+1}\mid_{\overline{ M}})$ satisfying the symmetry condition 
\begin{equation}
\label{eq.sym}
\forall \varrho_{0} \in Z, \forall \widehat {\varrho_{0}}, \widetilde{ \varrho_{0}}\in j^{-1}(\varrho_{0})\cap \text{Car} (\widetilde P), f(\widehat {\varrho_{0}}) = f(\widetilde{ \varrho_{0}})\end{equation}
then
$$\frac d {ds} f(j(\gamma(s))\mid_{s= s_{0}}= H_{\widetilde p}\mid_{j^{-1}( \gamma(s_{0}))} f(j^{-1}(\gamma( s_{0}))) $$
\end{enumerate}
\end{definition}
It is proved in~\cite{MeSj82} that under the assumption of no infinite order contact, through every point $\varrho_{o}\in {^b}T^* M \setminus \{0\}$ there exists a unique generalized bicharacteristic (which is furthermore a limit of bicharacteristics having only hyperbolic contacts with the boundary). This defines the flow $\Phi$. Finally remark that since $\widetilde {p}= p -\tau$ We have consequently a natural flow, $\varphi$ on $\text{Char}  {P} \subset {^b} T^*\Omega$ (the generalized flow of $p(z, \zeta)$) given by  
\begin{equation}
\Phi_{s}(t,\tau,z, \zeta)= (t-s, \tau, \varphi_{s}(z, \zeta))
\end{equation} 
\subsection{Wigner measures}
Consider functions
$a=a_i+a_\partial $ with $a_{i}\in C^\infty_{0}(T^* M)$, and $a_{\partial}\in C^\infty_{0}({\mathbb R}^{2d-1})$. Such symbols are quantized in the following way:
take $\varphi_{i}\in C^\infty_{0}(M)$  (resp $\varphi_{\partial}\in C^\infty_{0}({\mathbb R}^d)$) equal to $1$ near the $x$-projection of $\text{supp}(a_{i})$ (resp the $x$-projection of $\text{supp}(a_{\partial})$) and define
\begin{multline}
\label{eq3.1bis}
\text{Op}_{\varphi_{i}, \varphi_{\partial}}(a)(x, hD_{x})f= \frac 1 { (2\pi h)^d} \int e^{ i(x-y)\cdot \xi/h}a_{i}(x, \xi)\varphi_{i}(y)f(y) dy d\xi \\
+ \frac 1 { (2\pi h)^{ d-1}} \int e^{ i(x'-y')\cdot \xi'/h}a_{\delta}(x_{n}, x', \xi)\varphi_{\delta}(x_{n}, y')f(x_{n}, y') dy' d\xi'.   
\end{multline}
Remark that according to the symbolic semi-classical calculus, the operator $\text{Op}_{\varphi_{i}, \varphi_{\partial}}(a)$ does not depend on the choice of functions $\varphi_{i}, \varphi_{\partial}$, modulo operators on ${\mathcal {L}}^2$ of norms bounded by $O(h^\infty)$. As in the previous section, we shall in the sequel drop the index $\varphi_{i}, \varphi_{\partial}$\par
Denote by ${\mathcal A}$ the space of the operators which are a finite sum of operators obtained as above in suitable coordinate systems near the boundary and for $A\in {\mathcal {A}}$, by $a=\sigma(A)$ the semiclassical symbol of the operator $A$. For such functions $a$ we can define  $\kappa(a)\in C^0(Z)$ by
\begin{equation} 
\kappa(a)(\rho )=a(j^{-1}(\rho )) \label{eq:9}
 \end{equation} 
(the value is independent of the choice of $j^{ -1}( \rho)$ since the operator is tangential). \par
The set 
\begin{equation} 
\{\kappa(a),a=\sigma (A), A\in {\mathcal A}\} \label{eq:10}
 \end{equation} is a locally dense subset of
$C^0_{c}(Z)$.\par
\subsection{Elliptic regularity}
Consider a sequence $(f_{k})$ bounded in ${\mathcal {L}}^2= L^2_{\text{loc}}( {\mathbb R}_{t}; L^2( \Omega))$, solution of the equation (with $\lim_{k\rightarrow + \infty} h_{k}=0$)
\begin{equation}\left\{
\begin{aligned}
(ih_{k}\partial _{t}+ h_{k}^2 P) f_{k}&= o(h_{n})_{L^2_{\text{loc}}( {\mathbb R}_{t}; L^2( \Omega))}\\
u\mid_{\partial \Omega}&=0.
\end{aligned}\right.
\end{equation}
The same argument as in section~\ref{se.2.2.1} shows:
\begin{proposition}
If $a_{i}$ is equal to $0$ near $\text{Car}\widetilde{P}$ then 
\begin{equation}
\label{eq3.20}\lim_{k\rightarrow + \infty}\left(\text{Op}(a_{i})(x, h_{k}D_{x})f_{k}, f_{k}\right)_{ L^2}=0,
\end{equation}
\end{proposition}
and the analysis of the boundary value problem shows:
\begin{proposition}
If   $a_{\partial}$ is equal to $0$ near $Z$ (i.e. $ a_{i}$ is supported in the elliptic region)
 then
\begin{equation}
\label{eq3.21bis}\lim_{k\rightarrow + \infty}\left(\text{Op}(a_\partial)(x_{n}, x', h_{k}D_{x'})f_{k}, f_{k}\right)_{ L^2}=0.
\end{equation} 
\end{proposition}
\subsection{Definition of the measure}
The analog of proposition~\ref{prop2.1} is:
\begin{proposition}\label{prop:2.4} There exists a subsequence $(k_{p})$ and a Radon positive measure $\mu $ on $Z$  such that 
\begin{equation} 
\forall  Q\in{\mathcal A}\quad \lim_{p\rightarrow \infty }(Qf_{k_{p}},f_{k_{p}})_{L^2}=\langle\mu ,\kappa(\sigma(Q))\rangle.
 \end{equation}
\end{proposition}
The proof of this result in the interior of $\Omega$ is the same as in Section~2 and near a boundary point, it relies on the G\aa rding inequality for tangential operators (see G. Lebeau~\cite{Le96} for a proof in the classical context and~\cite{GeLe93-1, Bu97-1} for the semi-classical construction). As before, we denote again by $(f_{k})$ the extracted sequence.\par
\begin{proposition}[First properties of the measure $\mu$]
\label{prop3.2}
\begin{equation}
\label{eq3.27}\mu({\mathcal {H}})=0,
\end{equation}
\begin{equation}
\label{eq3.27'}\limsup_{k\rightarrow + \infty}| \left(\text{ Op}(a) {h_{k}} D _{x_{n}} f_{k}, f_{k}\right)_{ L^2}|\leq C \sup_{\varrho\in \text{ supp}( a)}|r|^{ 1/2} | a| .
\end{equation}
\end{proposition}
The relation~\eqref{eq3.27} is a simple consequence of the micro-local analysis of the boundary problem near a point $\varrho_{0}\in {\mathcal {H}}$, for which a parametrix for the solution can be written in terms of a semi-classical Fourier integral operator, by geometric optics methods. To prove~\eqref{eq3.27'} compute (with $\varphi \in C^\infty_{0}$ equal to $1$ near the $t$-projection of the support of $a$)
\begin{equation}
\label{eq3.40}\begin{aligned}\ 
\Big| \Big(\text{ Op}(a) {h_{k}}D _{x_{n}} f_{k}, f_{k}\Big)_{ L^2}\Big|
&\leq \| \text{ Op}(a)  {h_{k}} D _{x_{n}} f_{k}\|_{ L^2} \|\varphi(t) f_{k}\|_{ L^2}\\
&\leq  \bigl( {h_{k}} D _{x_{n}}\text{ Op}(a)^* \text{ Op}(a) {h_{k}}D _{x_{n}} f_{k}, f_{k} \bigr)_{ L^2}^{1/2}\|\varphi(t) f_{k}\|_{ L^2}\\
&\leq  \left(\text{ Op}(a)^* \text{ Op}(a) h_{k}^2D _{x_{n}}^2 f_{k}, f_{k} \right)_{ L^2}^{1/2}\|\varphi(t) f_{k}\|_{ L^2}+ o(1)\\
&\leq   \left(\text{ Op}(a)^* \text{ Op}(a) ( R-\widetilde{P})  f_{k}, f_{k} \right)_{ L^2}^{1/2}\|\varphi(t) f_{k}\|_{ L^2}+ o(1)\\
&\leq \left(\text{ Op}(a)^* \text{ Op}(a)  R  f_{k}, f_{k} \right)_{ L^2}^{1/2}\|\varphi(t) f_{k}\|_{ L^2}+ o(1) 
\end{aligned}
\end{equation}
and we obtain 
\begin{equation}
\label{eq3.41}\limsup_{k\rightarrow + \infty}\Big| \Big(\text{ Op}(a) {h_{k}}D _{x_{n}} f_{k}, f_{k}\Big)_{ L^2}\Big|
\leq C| \langle \mu, a^2 r \rangle | ^{ 1/2}\leq C \sup_{\varrho\in \text{ supp}( a)}|a||r|^{1/2}.
\end{equation}   
\subsection{Invariance of the measure}
Consider now a sequence $(f_{k})$ bounded in $L^2_{\text{loc}}( {\mathbb R}_{t}; L^2( \Omega))$, solution of the equation (with $\lim_{k\rightarrow + \infty} h_{k}=0$)
\begin{equation}\left\{
\begin{aligned}
(ih_{k}\partial _{t}+ h_{k}^2 P) f_{k}&= o(h_{k})_{L^2_{\text{loc}}( {\mathbb R}_{t}; L^2( \Omega))}\\
u\mid_{\partial \Omega}&=0
\end{aligned}\right.
\end{equation}
\begin{proposition}
\label{prop3.1}
Consider $q\in  C^\infty(T^*{\mathbb R}^{d+1}\mid_{\overline{ M}})$ satisfying the symmetry condition~\eqref{eq.sym}. In general $\{\widetilde p, q\}=-2 \xi_{n} \partial _{x_{n}}q+ \{r, q\}$ is not a function defined on $Z$ (because of the $\xi_{n}$ dependence). To obtain a function on $Z$, we take the convention 
\begin{equation}
\label{eq3.30bis}\{\widetilde p, q\}_{~}\overset{ def}=-2 \xi_{n} \partial _{x_{n}}q 1_{\varrho\notin {\mathcal {H}}}+ \{r, q\}.
\end{equation}
This function is $\mu$-integrable and, thanks to~\eqref{eq3.27}, $\mu$-almost everywhere continuous.

Then, with this convention, the measure $\mu$ satisfies
\begin{equation}
\label{eq3.26}\langle \mu, \{\widetilde p, q\}_{~} \rangle =0
\end{equation}
\end{proposition}

The proof of Proposition~\ref{prop3.1} is simply integration by parts (and some carefull study of the terms arising). We give it below: \par
Since in $M$, the equation~\eqref{eq3.26} is simply~\eqref{eq2.14}, we restrict the study to the case where $q$ is supported near a point $\varrho_{0}\in T^* \partial M$. Suppose first only that $q\in  C^\infty(T^*{\mathbb R}^{d+1}\mid_{\overline{ M}})$. From Malgrange preparation theorem, there exist functions $q_{0}(x_{n},x',\xi'), q_{1}(x_{n},x',\xi')\in C^\infty $ such that 
\begin{equation}
\label{eq3.30'}q\mid_{\text{Car}\widetilde{P}}= q_{0}\mid_{\text{Car} \widetilde{P}} + \xi_{n} q_{1}\mid_{\text{Car} \widetilde{P}}.
\end{equation}
Let $Q= \text{Op}(q_{0}) +\text{Op}(q_{1}) hD _{x_{n}}$ and compute ($P^*= P$)
\begin{equation}
 h_{k}^{ -1}\left((\widetilde{P}  Q - Q \widetilde{P})f_{k}, f_{k}\right)_{ L^2}.
\end{equation}
 Two integrations by part, \eqref{3.22} and the boundary condition $f_{k}\mid _{x_{n}=0}=0$ show that 
\begin{equation}
\label{eq3.29'}
\begin{aligned} h_{k}^{ -1}\left([\widetilde{P}, Q]f_{k}, f_{k}\right)_{ L^2}&= h_{k}^{ -1}\left(\widetilde{P}  Q f_{k}, f_{k}\right)_{ L^2}+ o(1)\\
&=  -i\left( Q_{1}\mid_{x_{n}=0} h_{k}D _{x_{n}}f_{k}\mid_{x_{n}=0}, h_{k}D _{x_{n}}f_{k}\mid_{x_{n}=0}\right)_{ L^2({\mathbb R}^{d-1}_{x'})}.
\end{aligned}
\end{equation}
On the other hand $[\widetilde {P}, Q]$ can be written under the form
$$\frac i{h_{k}}[ \widetilde {P}, Q]= A_{0}+ A_{1} h_{k}D_{x_{n}}+ A_{2}\widetilde{P} + hA_{3}$$
where $A_{0}$, $A_{1}$ and $ A_{2}$ are tangential operators, $A_{3}$ is differential of order at most $1$ in $D_{x_{n}}$, and on $\text{Car} \widetilde{P}$ we have 
\begin{equation}
\label{eq3.28}a_{0}+ a_{1} \xi_{n} = \{\widetilde p, q\}.
\end{equation}
From~\eqref{eq3.27} we deduce that $\mu$-almost everywhere 
\begin{equation}
\label{eq3.28'}a_{0}+ a_{1} \xi_{n} 1_{x_{n}>0} = \{\widetilde p, q\}.
\end{equation}
Consequently
\begin{multline}
\label{eq3.29}
h_{k}^{ -1}\left([\widetilde{P},  Q]f_{k}, f_{k}\right)_{ L^2}= h_{k}^{ -1}\left(\widetilde{P}  Qf_{k}, f_{k}\right)_{ L^2}+ o(1)\\
=\left(-i (A_{0} + A_{1} {h_{k}}D _{x_{n}} + A_{2}\widetilde{P}+ o(1)) f_{k}, f_{k}\right)_{ L^2}+ o(1).
\end{multline}
Passing to the limit in~\eqref{eq3.29} we obtain
\begin{equation}
\label{eq3.29''}
 \lim_{k\rightarrow + \infty}\left( (A_{0} + A_{1} {h_{k}}D _{x_{n}} + A_{2}(\widetilde{P}) f_{k}, f_{k}\right)_{ L^2}\\
= \langle \mu, a_{0}\rangle + \lim_{k\rightarrow + \infty}\left( A_{1} {h_{k}}D _{x_{n}} f_{k}, f_{k}\right)_{ L^2}. 
\end{equation} 
Take $\varepsilon>0$ and $\varphi\in C^\infty_{0}(]-1, 1[)$ equal to $1$ near $0$. Decompose
\begin{multline}
\label{eq3.32bis} A_{1}=  (1-\varphi(\frac {x_{n}} \varepsilon ))A_{1}
+\text{Op}\left(\varphi (\frac {x_{n}} \varepsilon )\varphi(\frac { r(x_{n},x',\xi')} { \varepsilon})\right)A_{1} \\
+\text{Op}\left(\varphi(\frac {x_{n}} \varepsilon )( 1-\varphi(\frac { r(x_{n},x',\xi')} { \varepsilon}))\right)A_{1}.  
\end{multline}
The first term in the right hand side of~\eqref{eq3.32bis} is supported in the interior of $\Omega$; its contribution to the limit in~\eqref{eq3.29''} is equal to 
\begin{equation}
\label{eq3.33bis}\langle \mu, (1- \varphi(\frac {x_{n}} \varepsilon))a_{1}\xi_{n}\rangle 
\end{equation}
The contribution of the second term is, according to~\eqref{eq3.27'}, smaller than
\begin{equation}
\label{eq3.35}C\sup_{\varrho\in \text{supp}(\varphi(\frac {x_{n}} \varepsilon )(\varphi(\frac { r({x_{n}},x',\xi')} { \varepsilon}))) } | r|^{1/2} | a_{1}|\leq C \varepsilon^{ 1/2}
\end{equation}
and the contribution of the last term is smaller than
\begin{multline}
\label{eq3.35'}\| h_{k}D _{{x_{n}}}f_{k}\| \Big\|A_{1}^*\text{Op}\left(\varphi(\frac {x_{n}} \varepsilon )( 1-\varphi(\frac { r({x_{n}},x',\xi')} { \varepsilon}))\right)^* f_{k}\Big\| \\
\leq C \langle \mu, |a_{1}|^2\varphi^2(\frac {x_{n}} \varepsilon )( 1-\varphi(\frac { r({x_{n}},x',\xi')} { \varepsilon}))^2\rangle ^{ 1/2} + o(1) 
\end{multline}
Passing to the limit $\varepsilon\rightarrow 0$ we obtain that the contribution of the first term is equal to 
\begin{equation}
\label{eq3.36}\langle \mu , a_{1}\xi_{n} 1_{{x_{n}}>0}\rangle,
\end{equation}
the contribution of the second term is (according to~\eqref{eq3.35}) equal to $0$ 
and the contribution of the last term is, according to~\eqref{eq3.27}, smaller than
\begin{equation}
\label{eq3.36'}\langle \mu, a_{1}^2 1_{{x_{n}}=0}1_{ r\neq 0}\rangle =\langle \mu , a_{1}^2 1_{\varrho\in {\mathcal {H}}}\rangle=0.
\end{equation}
Finally we have proved 
\begin{equation}
\label{eq3.50}\lim_{n\rightarrow + \infty}-i\left( Q_{1}\mid_{{x_{n}}=0} h_{k}D _{{x_{n}}}f_{k}\mid_{{x_{n}}=0}, h_{k}D_{{x_{n}}}f_{k}\mid_{{x_{n}}=0}\right)_{ L^2({\mathbb R}^{d-1}_{x'})}
=-i\langle \mu, \{\widetilde{p}, q\}\rangle. 
\end{equation} 
But, if $q$ satisfies the symmetry condition~\eqref{eq.sym}, the function $q\mid_{{x_{n}}=0}$ is independent of $\xi_{n}$ on $\text{Car}P$. Hence $q_{1}\mid_{{x_{n}}=0}=0$ on ${\mathcal {H}}$ and consequently on $\overline{H}= {\mathcal {H}} \cup {\mathcal {G}}$; and the left hand side in~\eqref{eq3.50} tends to $0$, which proves Proposition~\ref{prop3.1}.
\begin{proposition}[see~\cite{BuGe96} and~\cite{Mi97}]
\label{prop3.3}We have 
\begin{equation}
\label{eq3.52} \mu({\mathcal {G}}^{2,+})=0
\end{equation}
\end{proposition}
Consider a point $\varrho_{0}\in {\mathcal {G}}^{2,+}$. Apply~\eqref{eq3.50} to a family of functions $q= \xi_{n} \times q_{\varepsilon}$ with
\begin{equation}
q_{\varepsilon}=\varphi(\frac {x_{n}}{\varepsilon^{ 1/3}}) \varphi(\frac{(r({x_{n}},x',\xi'))}{\varepsilon})a({x_{n}},x',\xi').
\end{equation}
 Then we get 
\begin{multline}
\label{eq3.53}\textstyle \lim_{n\rightarrow + \infty}\left( \varphi(\frac{(r(0,x',h_{k}D_{x'}))}{\varepsilon})a(0,x',h_{k}D_{x'}) h_{k}D _{{x_{n}}}f_{k}\mid_{{x_{n}}=0}, h_{k}D _{{x_{n}}}f_{k}\mid_{{x_{n}}=0}\right)_{ L^2({\mathbb R}^{d-1}_{x'})}\\
=\langle \mu,-2 \xi_{n}^2 \partial _{{x_{n}}}(\varphi(\frac {x_{n}}{\varepsilon^{ 1/3}}) \varphi(\frac{(r({x_{n}},x',\xi'))}{\varepsilon})a({x_{n}},x',\xi'))\rangle\\
- \langle \mu, \partial _{{x_{n}}}r \varphi(\frac{x_{n}}{\varepsilon^{ 1/3}}) \varphi(\frac{(r({x_{n}},x',\xi'))}{\varepsilon})a({x_{n}},x',\xi')+ \xi_{n} \{r, q_{\varepsilon}\}'  \rangle, 
\end{multline}
where $\{r,q_{\varepsilon}\}'$ is the Poisson bracket with respect to the $x', \xi' $ variables. On the support of the measure $\mu$, $\xi_{n}^2= r({x_{n}},x',\xi')$. Hence we can apply the dominated convergence theorem and obtain that the right hand side in~\eqref{eq3.53} tends to 
\begin{equation}
\label{eq3.54} \langle \mu, -\partial _{{x_{n}}}r(0,x',\xi')a({x_{n}},x',\xi')1_{{x_{n}}=0}1_{r=0}\rangle=\langle \mu, -\partial _{{x_{n}}}r(0,x',\xi')a({x_{n}},x',\xi')1_{\rho \in {\mathcal {G}}}\rangle 
\end{equation}
According to the assumption $\varrho_{0}\in {\mathcal {G}}^{2,+}$,  $\partial _{{x_{n}}}r >0$ at the point $\varrho_{0}$ . If the support of $a$ is chosen small enough so that $\partial _{{x_{n}}}r >0$ on this support, then the right hand side in~\eqref{eq3.54} is non positive. 
On the other hand by G\aa rding inequality the limit on the left hand side is non negative. Both sides are then equal to $0$. This implies Proposition~\ref{prop3.3}.\par
It is now possible to prove as in ~\cite[Th\protect{\'e}or\protect{\`e}me~1]{BuLe01} (see also~\cite{Bu02}), by measure theory methods, that the invariance of the measure $\mu$ along the generalized bicharacteristic flow is equivalent to Propositions~\ref{prop3.1} and~\ref{prop3.3} (in fact the proof of this result is presented in~\cite[Section 3.3]{BuLe01} for classical measures, in the more general context of systems, but the proof for semi-classical measures is the same word by word).

\subsection{Proof of Theorem~\protect\ref{th2}}
All that remains to do to complete the proof of Theorem~\ref{th2} in the case of a Dirichlet boundary value problem is to construct a sequence of initial data $(u_{0,n})$ and a sequence $h_{n}; \lim_{n\rightarrow + \infty} h_{n}=0$, such that the sequence of solutions of the semi-classical Schr\"odinger equations admits  
\begin{equation}\label{eqmespr}
dt \otimes\delta_{\tau= p(z_{0}, \zeta_{0})}\otimes \delta_{(z, \zeta)=(\varphi_{-t}(z_{0}, \zeta_{0}))}
\end{equation}
 as Wigner measure. \par
In the case where the bicharacteristic starting from $(t_{0}=0, \tau_{0}= p(z_{0}, \zeta_{0}), z_{0}, \zeta_{0})$ has an interior point $(t_{1}, \tau = \tau_{0}, z_{1}\in \Omega, \zeta_{1})$, we perform the construction as in the previous section, since by finite speed of propagation (modulo ${\mathcal {O}}( h^ \infty)$), the boundary is not seen,~\eqref{eqmespr} is satisfied close to $(t_{1}, \tau=\tau_{0}, z_{1}, \zeta_{1})$. Using the propagation result, we deduce that~\eqref{eqmespr} is satisfied everywhere.\par
 In the case where the bicharacteristic starting from $(t_{0}=0, \tau=\tau_{0}, z_{0}, \zeta_{0})$ has no interior point, we know that it can be approximated by bicharacteristics $\gamma_{k}$ which have an interior point (see~\cite{MeSj78, MeSj82}). For these  bicharacteristics, we can construct sequences of initial data $(u_{n,k})$ associated to $(h_{n,k}); \lim_{n\rightarrow + \infty} h_{n,k}=0$. Taking $(u_{n_{k},k})$ with $n_{k}$ large enough, as initial data matches our aim.\par
 
  The rest of the proof of the estimate~\eqref{eq3.1} in the case of a boundary value problem is now the same as in Section~\ref{sec.3}.

\section{Smoothing effect}\label{sec.2}
In this section we prove a weaker smoothing effect for a class of trapping obstacles.
\subsection{Assumptions}
Consider $\Theta\subset {\mathbb R}^d$ a compact smooth obstacle whose complement, $\Omega = \Theta^ c$ is connected. Let  $\Delta_{D}$ be the Laplace operator acting on $L^2(\Omega)$, with domain $D= H^2(\Omega)\cap H^1_{0}(\Omega)$. Denote, for $u_{0}\in L^2(\Omega)$, by $e^{-it\Delta_{D}}u_{0}=u$ the solution of the Schr\"odinger equation with Dirichlet boundary conditions:
\begin{equation}
\label{eqLS}\left\{\begin{aligned} (i\partial_{t}- \Delta)u &=0 \text { on }{\mathbb R}_{t}\times \Omega,\\
u\mid_{\partial \Omega}&=0,\\
u\mid _{t=0} &= u_{0}.
\end{aligned} \right.
\end{equation}
We suppose that $\Theta= \cup _{i=1}^N \Theta_{i}\subset {\mathbb R}^d$ is the union of a finite number of strictly convex obstacles, $\Theta_{i}$ satisfying:
\begin{itemize}
\item For any $1\leq i,j,k\leq N$, $i\neq j$, $j\neq k$, $k\neq i$, one has 
\begin{equation}
\label{eq4}{\text{Convex Hull}}(\Theta_{i}\cup \Theta_{j}) \cap \Theta_{k}= \emptyset.
\end{equation}
\item Denote by $\kappa$ the infimum of the principal curvatures of the boundaries of the obstacles $\Theta_{i}$, and $L$ the infimum of the distances between two obstacles. Then if $N>2$ we assume that $\kappa L >N$ (no assumption if $N=2$). 
   \end{itemize}
\begin{remarque}
\label{re2.30} If there are only two obstacles, then the assumptions are automatically fulfilled. The first assumption is essentially technical, whereas the second one is an assumption about the strong hyperbolicity of the dynamical system given by the billiard flow.
\end{remarque}
 In this case, since there are trapped trajectories (for example any line minimizing the distance between two obstacles is trapped), we have shown in Section~\ref{sec.4}  that the plain smoothing effect $H^{1/2}$ does not hold. However, the result below (a more precise version of Theorem~\ref{th6}) shows that the smoothing effect with a logarithmic loss still holds.
 \begin{theoreme}
\label{th4}Under the assumptions above, for any $\chi \in C^\infty _{0}({\mathbb R}^d)$ there exists $C>0$ such that the solution of 
\begin{equation}
\label{eq2.100}\left\{
\begin{aligned} (i\partial_{t}- \Delta)u &=0 \text { on }{\mathbb R}_{t}\times \Omega,\\
u\mid _{t=0} &= u_{0},\\
u\mid_{\partial \Omega}&=0
\end{aligned}\right.
\end{equation}
and the solution of
\begin{equation}
\label{eq2.100bis}\left\{
\begin{aligned} (i\partial_{t}- \Delta)v &=\chi f; \text { $\chi f$ compactly supported in time},\\
v\mid _{t<<0} &=0,\\
 v\mid_{\partial \Omega}&=0
\end{aligned}\right. 
\end{equation}
satisfy:
\begin{equation}
\label{eqsmoothbis}\begin{aligned}
\| \chi u \|_{L^2({\mathbb R}_{t}; H_{D}^{1/2,- }(\Omega))}&\leq C\|u_{0}\|_{L^2( \Omega)},\\
\| \chi v \|_{L^2({\mathbb R}_{t}; H_{D}^{1/2,-}(\Omega))}&\leq C \|\chi f\| _{L^2({\mathbb R}_{t}; H_{D}^{-1/2,+}(\Omega))},
\end{aligned}
\end{equation}
where $H_{D}^{1/2,-}= D((\text{Id}- \Delta_{D})^{1/4}\log^{ -1/2}(2\text{Id}- \Delta_{D}))$ and $H^{-1/2, +}=(H_{D}^{1/2,-})'$. In particular
$$\forall \varepsilon>0, H^{-1/2+\varepsilon}(\Omega)\subset H_{D}^{-1/2,+}\subset H^{-1/2}(\Omega), \qquad H^{1/2}(\Omega)\subset H_{D}^{1/2,-}\subset H^{1/2- \varepsilon}(\Omega)$$
with continuous injections.
\end{theoreme}
\begin{remarque}\label{re4.7}
In the case where there exist an elliptic (stable) periodic trajectory, it is possible to construct quasi modes with compact support, i.e. functions $(e_{n})_{n\in {\mathbb N}}$ with compact supports associated to a particular sequence $(\lambda_{n})\rightarrow + \infty$ and satisfying
\begin{equation}\begin{aligned}
\label{eq2.100ter}-\Delta e_{n}&= \lambda_{n}e_{n}+ r_{n},\\
\|r_{n}\|_{{H^N}}&\leq C_{N,M} \lambda_{n}^{-M}, \ \forall N, M \in {\mathbb N}.
\end{aligned}
\end{equation}
From this we deduce easily that the sequence of solutions of the Schr\"odinger equation with initial data $(e_{n})$ is, for any $\varepsilon>0, s>0$, not bounded in $L^1\left([0,\varepsilon[; H^s_{{\text{loc}}}\right)$;  which implies that no smoothing effect at all is true any more. Under the assumptions of Theorem~\ref{th4}, the periodic trajectories are hyperbolic (unstable), which forbids the construction of such well localized quasi-modes.
\end{remarque}
Theorem~\ref{th4} is deduced  from the following estimate of the cut-off resolvent:
\begin{proposition}
\label{prop4}
Suppose that the obstacle $\Theta$ satisfies the assumptions in Theorem~\ref{th4} above. Then the resolvent of the operator $\Delta_{D}$, $(-\Delta_{D}- \lambda)^{-1}$ (which is analytic in ${\mathbb C}\setminus {\mathbb R}^+$) satisfies:
\begin{equation}
\label{eq2.15}\begin{gathered}\forall \chi \in C^ \infty_{0}({\mathbb R}^2), \exists C>0; \forall \lambda \in {\mathbb R}, 0<\varepsilon<<1,\\
\| \chi ( -\Delta_{D}- (\lambda\pm i \varepsilon))^{-1}\chi\|_{L^2\rightarrow L^2}\leq \frac {C \log(2+ |\lambda|)} { 1 + \sqrt{ |\lambda| }}.
\end{gathered}
\end{equation}
\end{proposition}
We are going to prove  this estimate for $\lambda>>1$. The proof for $|\lambda|<<1$ can be found in~\cite[Annexe B.2]{Bu98}; whereas the result for $c\leq |\lambda|\leq C$ follows from the Rellich uniqueness Theorem (see~\cite{LaPh89}or \cite[Annexe B.1]{Bu98}) and the result for $\lambda <-\varepsilon$ is clear because in this case the operator is semi-classically elliptic.
\par Let us perform a change of variables $\lambda = \tau^2$ and consider $ \chi ( -\Delta_{D}- (\tau^2))^{-1}\chi$ which is holomorphic in $\{\Im \tau <0\}$ and satisfies there (according to the standard estimate for self adjoint operators),
\begin{equation}
\label{eq2.10}\|( -\Delta_{D}- (\tau^2))^{-1}\|_{L^2(\Omega)\rightarrow L^2(\Omega)}\leq \frac 1 { |\tau| |\Im \tau|}.
\end{equation}
M. Ikawa proved in~\cite{Ik82, Ik83} and more precisely in~\cite[Theorem 2.1]{Ik88} (see also the work by C. G\'erard~\cite{Ge88} where such an estimate is implicit) that under the assumptions above, the following estimate on the cut-off resolvent holds:
\begin{theoreme}[Ikawa, \protect{~\cite[Theorem 2.1]{Ik88}}]
\label{thIk}
The cut-off resolvent: $ \chi ( -\Delta_{D}- (\lambda\pm i \varepsilon))^{-1}\chi$
 admits a holomorphic continuation in a strip of the upper half plane 
\begin{equation}
\label{eq2.11}\{ \tau\in {\mathbb C}; |\tau| >1, \Im \tau \leq \alpha\}, \alpha>0
\end{equation}
and satisfies there (for a large $N$):
 \begin{equation}
\label{eq2.12}\|\chi ( -\Delta_{D}- (\tau^2))^{-1}\chi \|_{L^2(\Omega)\rightarrow L^2(\Omega)}\leq C |\tau| ^N
\end{equation}
\end{theoreme}
\begin{remarque}
In~\cite[Theorem 2.1]{Ik88} the proof is done with the additional assumption that the dimension of space is equal to $3$ (which is the relevant dimension the author had in mind for applications to the wave equation). However the proof could be equally performed in any space dimension $d\geq 2$ (see~\cite{Ge88} in the case $N=2$, $d\geq 2$).
\end{remarque}
Using~\eqref{eq2.10} and~\eqref{eq2.12} (and writing $\tau=h^{-1}z$, $z\sim 1$, $h\rightarrow 0$), one easily sees, with 
\begin{equation}
f(h,z)= \left( \chi ( -h^2\Delta_{D}- z^2)^{-1} \chi u,v\right) _{{L^2(\Omega)}}; u, v \in L^2( \Omega)
\end{equation}
 that~\eqref{eq2.15} for large $|\lambda|$ follows from Theorem~\ref{thIk} and the following  semi-classical maximum principle (a variant of Phr\"agmen Lindel\"of principle) adapted from the work by Tang--Zworski~\cite{TaZw00}:
\begin{lemme}
\label{le2.1} Suppose that $f(h,z)$ is a family of holomorphic functions defined for $0<h<1$ in a neighbourhood of 
\begin{equation}
\Omega(h)= [1/2, 3/2]\times i[h\alpha, -h\alpha]
\end{equation}
such that 
\begin{equation}\label{eq2.20}
  \begin{aligned}|f(h,z)|&\leq Ch^{-M} \text{on $\Omega(h)$},\\
|f(h,z)|&\leq \frac 1 {|\Im z|} \text{on $\Omega(h)\cap \{\Im z<0\}$}.
\end{aligned}
\end{equation}
Then there exists $h_{0}>0$, $C>0$, such that for any $0<h<h_{0}$
\begin{equation}|f(h,z)|\leq C \frac{\log (h^{-1})} h \text{on $[4/5, 6/5]$}
\end{equation}
\end{lemme}
To prove this lemma, first consider the function
\begin{equation}
\varphi(z,h)= (\pi h)^{-1/2}\int e^{-\frac{(x-z)^2} h}\Psi(x) dx
\end{equation}
where $\Psi \in C^\infty_{0}(]2/3, 4/3[)$ is non negative and equal to $1$ in $[3/4, 5/4]$. Then the function $\varphi(z,h)$ satisfies:
\begin{enumerate}
\item $\varphi (z,h)$ is holomorphic in $\Omega(h)$,
\item $|\varphi(z,h)|\leq C$ in $\Omega(h)$,
\item $|\varphi(z,h)|\geq c>0$ in $[4/5, 6/5]$,
\item $|\varphi(z,h)|\leq Ce^{-c/h}$ on $\Omega(h)\cap \{ |\Re z- 1| \geq 1/2\}$.   
\end{enumerate}
Then apply the maximum principle to the function $g(z,h)= e^{-iN\log(h) z/h}\varphi(z,h) f(z,h)$ on the domain 
\begin{equation}
\widetilde { \Omega}(h)= [1/2, 3/2]\times i[h\alpha, -h/\log(h^{-1})].
\end{equation}
Using the bounds~\eqref{eq2.20} on $f$ and the properties of $\varphi$ above, we can estimate $g$ by
\begin{equation}
\begin{aligned}
|g(z,h)|&\leq C h^{N\alpha -M} \text{on $\partial \widetilde{\Omega}(h)\cap \{ \Im z = h\alpha \}$ },\\
|g(z,h)|&\leq C_{N} e^{-c/h} \text{ on $\partial \widetilde{ \Omega}(h)\cap \{ \Re x \in \{1/2, 3/2\}\}$},\\
|g(z,h)|&\leq C_{N} \frac{\log (h^{-1})} h \text{on $\partial \widetilde{\Omega}(h)\cap \{ \Im z = \frac{-h}{ \log(h^{-1})}\}$}.
\end{aligned}
\end{equation}
Taking $N$ large enough and applying the maximum principle we get 
\begin{equation}
|g(z,h)|\leq C' \frac{\log (h^{-1})} h \text{ on $\widetilde{\Omega}(h)$ }
\end{equation}
which implies 
\begin{equation}
|f(z,h)|\leq C' \frac{\log (h^{-1})} h \text{ on $[4/5, 6/5]$, }
\end{equation}
and ends the proof of Lemma~\ref{le2.1}.\par
We deduce from~\eqref{eq2.15}:
\begin{equation}
\label{eq2.110}
\| \chi ( -\Delta_{D}- (\lambda\pm i \varepsilon))^{-1}\chi\|_{H^{-1/2,+}\rightarrow H^{1/2,-}}\leq C.
\end{equation}
Indeed, for bounded $\lambda$, integrations by parts show that we can in fact replace $H^{1/2,-}$ by $H^1_{0}( \Omega)$ and $H^{-1/2,+}$ by $H^{-1}( \Omega)$ and for large $\lambda$ we decompose, with $\Psi \in C^\infty_{0}(]1/2,2[)$ equal to $1$ close to $1$,
\begin{equation}
u= (P-\lambda)^{-1}\chi f = \Psi(\frac {- \Delta_{D}} {\lambda}) u + (1- \Psi(\frac {- \Delta_{D}} {\lambda})) u.
\end{equation}
We get by the functional calculus of self adjoint operators 
\begin{equation}
\|(1- \Psi(\frac {- \Delta_{D}} {\lambda})) u\|_{H^1_{0}( \Omega)}\leq C \|(1- \Psi(\frac {- \Delta_{D}} {\lambda})) \chi f\|_{H^{-1}( \Omega)}.
\end{equation}
On the other hand, the function $ v= \Psi(\frac {- \Delta_{D}} {\lambda}) u$ satisfies
\begin{equation}\label{eq.trunc} (P- \lambda)v=\Psi(\frac {- \Delta_{D}} {\lambda}) \chi f.
\end{equation}
If $\widetilde{\chi}\in C^\infty_{0}( {\mathbb R}^d)$ is equal to one on the support of $\chi$, we have modulo negligible terms
\begin{equation}\label{eq.commute}
\widetilde \chi \Psi(\frac {- \Delta_{D}} {\lambda}) \chi= \Psi( \frac {- \Delta_{D}} {\lambda})\chi 
\end{equation}
because 
$$\widetilde \chi \Psi(\frac {- \Delta_{D}} {\lambda}) \chi- \Psi( \frac {- \Delta_{D}} {\lambda})\chi= [\widetilde \chi,\Psi(\frac {- \Delta_{D}} {\lambda})] \chi$$
and  on the support of $\nabla \widetilde \chi$, the operator $P$ is a differential operator and consequently $\Psi( \frac {- \Delta_{D}} {\lambda})$ is a pseudodifferential operator on this set (see for example Sect.4 of~\cite{Sj97}).\par
 According to~\eqref{eq.trunc}, ~\eqref{eq.commute} and Lemma \ref{le2.1}, we get
\begin{equation}
\left(\frac{\sqrt{ |\lambda|}}{\log(2 + |\lambda|)}\right)^{1/2}\|\chi \Psi(\frac {- \Delta_{D}} {\lambda}) u \|_{L^2} \leq C \left(\frac{\log(2 + |\lambda|)}{\sqrt{ |\lambda|}}\right)^{1/2}\|\Psi(\frac {- \Delta_{D}} {\lambda}) \chi f\|_{L^2} 
\end{equation}
to replace the weights in $\lambda$ above by the $H^{\pm 1/2, \mp}$-norms (i.e. to replace the weights in $\lambda$ by weights in $- \Delta_{D}$), it is enough to check that modulo negligible terms, if $\widetilde \Psi $ is a fonction equal to $1$ on the support of $\Psi$, 
\begin{equation}
\widetilde \Psi(\frac {- \Delta_{D}} {\lambda}) \chi \Psi( \frac {- \Delta_{D}} {\lambda})=\chi \Psi( \frac {- \Delta_{D}} {\lambda}) 
\end{equation}
which follows from the same arguments as above.\par

Following~\cite{BuGeTz02}, Theorem~\ref{th4} is now a consequence of~\eqref{eq2.110}. For the sake of completeness and since the argument is short, we recall it: firstly remark that by $TT^*$ argument it suffices to study the second (inhomogeneous) case. Indeed denote by $T= \chi e^{-it\Delta_{D}}$. The continuity of $T$ from $L^2$ to $L^2({\mathbb {R}}_{t}; H^{1/2,-})$ is equivalent to the continuity of the adjoint operator 
\begin{equation}
\label{eq2.5}T^*f= \int_{{\mathbb R}}e^{is\Delta_{D}} \chi f(s) ds
\end{equation} 
from $L^2({\mathbb {R}}_{t}; H^{-1/2,+})$ to  $L^2$, which in turns is equivalent to the continuity of the operator $T T^*$ from $L^2({\mathbb {R}}_{t}; H^{-1/2,+})$ to $L^2({\mathbb {R}}_{t}; H^{1/2,-})$. But
\begin{equation}
\label{eq2.6}\begin{aligned} T T^* f(t)&= \int_{{\mathbb R}}\chi e^{i(s-ts)\Delta_{D}} \chi f(s) ds\\
&=  \int_{s<t}\chi e^{i(s-t)\Delta_{D}} \chi f(s) ds + \int_{t<s}\chi e^{i(s-t)\Delta_{D}} \chi f(s) ds
\end{aligned}
\end{equation}
and (by time inversion), it clearly suffices to prove the continuity of any one of the terms in the right hand side, which is the second (inhomogeneous) part of Theorem~\ref{th4}.\par
Consider now $(v,f)$ solution of~\eqref{eq2.100bis}. By translation invariance we can suppose that $f$ (and hence $v$) is supported in $\{ t>0\}$. The Fourier transforms of $v$ and $f$ are (according to the support property) holomorphic in the set $\{\Im z <0\}$ and satisfy there, according to~\eqref{eq2.100bis}
\begin{equation}
\label{ea2.7}(-z+ \Delta) \hat u(z, \cdot)= \chi \hat f (z, \cdot).
\end{equation}
Taking $z= \tau -i\varepsilon$, $\tau\in {\mathbb R}$ and having $\varepsilon$ tend to $0$, using~\eqref{eq2.110}, we get
\begin{equation}
\label{eq2.8}\|\chi \hat u \|_{L^2({\mathbb R}_{\tau}; H^{1/2,-})}\leq C \|\chi \hat f \|_{L^2({\mathbb R}_{\tau}; H^{-1/2,+})}
\end{equation}
and since the Fourier transform is an isometry on $L^2({\mathbb R}; H)$ if $H$ is a Hilbert space, we get~\eqref{eqsmoothbis}.\par

Finally, as in~\cite{BuGeTz02}, we can deduce from Theorem~\ref{th4}:
\begin{theoreme}[Global existence for 2-d defocusing NLS]\label{th5}
Consider $\Theta\subset {\mathbb R}^2$ an obstacle which is  the union of $N$ strictly convex obstacles satisfying the assumptions above. Denote by $\Omega= \Theta^c$ its complement and
let $P$ be a polynomial with real coefficients. 
For every $u_{0}\in H_{0}^{1}(\Omega)$, there exists a unique maximal
solution $u\in C(I,H_{0}^{1}(\Omega))$ of the equation
\begin{eqnarray}\label{1.7}
i\partial_{t}u+\Delta u = P'(|u|^2)u,\quad u(0,x)=u_{0}(x).
\end{eqnarray}
Moreover we have:
\begin{enumerate}  
\item[(i) ]  
If $\|u_{0}\|_{H_{0}^{1}(\Omega)}$ is bounded from above, the length
of $I\cap \mathbb{R}_{\pm}$ 
is bounded from below by a positive constant.
\item[(ii) ]  
For any finite $p$, $u\in L^{p}_{loc}(I,L ^{\infty}(\Omega))$.
\item[(iii) ]   
If $P(r) \longrightarrow +\infty$ as $r \longrightarrow +\infty$, $I=\mathbb{R}$.
\item[(iv) ]  
If $u_{0}\in H_{D}^{s}(\Omega)$ for some $s>1$, 
$u\in C(I,H_{D} ^{s}(\Omega))$. 
In particular if $u_{0}\in C_{0}^{\infty}(\Omega)$,
$u\in C^{\infty}(I\times \Omega)$.
\end{enumerate}
\end{theoreme}

\def\cprime{$'$}

\end{document}